\documentclass[12pt]{amsart}
\usepackage{amsfonts,amssymb,latexsym,epsfig,amsmath}
\usepackage{amssymb}
\usepackage{amsmath}
\numberwithin{equation}{section}
\newtheorem{thm}{Theorem}[section]
\newtheorem{lemma}{Lemma}[section]

\newtheorem{remark}{Remark}[section]

\newcommand{\sn}{~\hfill\hbox{$\square$}}
\newcommand{\n}{\hbox{$\scriptstyle\hbox{\rm l\kern-0.22em N}$}}
\newcommand{\N}{\hbox{$ I\!\!N$}}

\renewcommand{\r}{{\mathbb R}}
\def\r{{\mathbb R}}
\newcommand{\R}{\r^{N}}

\renewcommand{\O}{\mathcal{O}}
\newcommand{\e}{\varepsilon}
\newcommand{\tr}{\hbox{\rm tr}}

\def\PS #1 #2{\langle #1, #2 \rangle}

\def\1{1\!\!1}

\def\dO{{\partial \O}}
\def\dP{{\partial P}}
\def\sym{S^N}

\def\Z{\mathbb Z}

\def\wb{\bar w}

\def\e{\varepsilon}

\def\Ex{I\!\!E_x}

\def\Ft{{\tilde F_\varepsilon}}
\def\gb{\bar g}
\def\ue{u^\varepsilon}
\def\Fe{F_\varepsilon}
\def\ve{v^\varepsilon}
\def\xe{x_\e}
\def\te{t_\e}

\def\xoe{\varepsilon^{-1}x}
\def\toe{\varepsilon^{-2}t}

\def\FF{\mathcal{F}}
\def\OO{\mathcal{O}}

 1

\begin{document}

\title[]{\bf
On homogenization problems for fully
nonlinear equations with oscillating Dirichlet
boundary conditions}

\author[]
{G. Barles$^{(1)}$
and
E. Mironescu$^{(2)}$}

\addtocounter{footnote}{1} \footnotetext{Laboratoire de
Math\'ematiques et Physique Th\'eorique (UMR~7350). F\'ed\'eration Denis Poisson (FR~2964) Universit\'e de Tours.
Facult\'e des Sciences et Techniques, Parc de Grandmont, 37200
Tours, France. E-mail address: barles@lmpt.univ-tours.fr }

\addtocounter{footnote}{1} \footnotetext{Institut Camille Jordan UMR CNRS 5208, Universit\'e de Lyon,
Ecole Centrale de Lyon,
36 avenue Guy de Collongue,
69134 Ecully cedex, France. elisabeth.mironescu@ec-lyon.fr
}

\date{}

 \begin{abstract}
We study two types of asymptotic problems whose common feature - and difficulty- is to exhibit oscillating Dirichlet boundary conditions : the main contribution of this article is to show how to recover the Dirichlet boundary condition for the limiting equation. These two types of problems are (i) periodic homogenization problems for fully nonlinear, second-order elliptic partial differential equations set in a half-space and (ii) parabolic problems with an oscillating in time Dirichlet boundary condition. In order to obtain the Dirichlet boundary condition for the limiting problem, the key step is a blow-up argument near the boundary which leads to the study of Dirichlet problems set on half-space type domains and of the asymptotic behavior of the solutions when the distance to the boundary tends to infinity.
\end{abstract}



\maketitle

\section{Introduction}

In this article, we are interested in two different types of problems : the first one is homogenization problems in half-space domains of the following form : find the behavior when $\e \to 0$ of the solutions $\ue$ of the Dirichlet problem
\begin{equation}\label{hom-P}
\left\{
\begin{array}{rcl}
F_\e (D^2 \ue, D \ue , \ue , \xoe,  x)  & = & 0  \quad \hbox{in  }P\, ,\\
\displaystyle u^\e  & = &g (\xoe , x)
\quad \hbox{on  }\dP \; .
\end{array}
\right.
\end{equation}
where $P:=\{x\in \R;\ x \cdot e > 0\}$ for some $e \in \R$, $|e|=1$, $D\ue$ and $D^2\ue$ denote respectively the
gradient and Hessian matrix of the real-valued solution $u^\e$. The functions $F_\e(M,p,t,y,x)$ and $g(y,x)$ are real-valued, continuous functions defined on  $\sym \times \R \times \r \times \R\times P$ and $\R\times \dP $ respectively, where $S^N$ is the space of $N\times N$ symmetric matrices; our first key assumption is that both $F_\e$ and $g$ are assumed to be periodic in $y$. Since our approach relies on the Maximum Principle and viscosity solutions' theory, the functions $F_\e$ are also assumed to satisfy some ellipticity requirement.

The second type of problems concerns parabolic equations with oscillating in time Dirichlet boundary conditions, namely
\begin{equation}\label{osc-evol}
 \left\{
 \begin{array}{rcl}
\ue_t + \Ft (D^2 \ue ,  D \ue, x,t)  & = & 0
\quad \hbox{in  }\O \times (0,T) \, ,\\
\noalign{\vskip6pt}
\ue (x,t) & = & g(x, \toe) \quad \hbox{on  }\dO \times (0,T) \; ,\\
\noalign{\vskip6pt}
\ue (x,0) & = & u_0(x) \quad \hbox{in  }\O
\end{array}
 \right.
\end{equation}
where $\O$ is a smooth, bounded domain, the nonlinearities $\Ft$ and the data $u_0,g$ are continuous functions, the functions $\Ft$ satisfying a suitable ellipticity condition, and $g(x,\tau)$ being $1$-periodic in $\tau$ for any $x\in \dO$.

The aim of this article is to identify the Dirichlet boundary condition $\gb$ for the limiting problems and we are going to concentrate on this issue, without providing in details the full asymptotic results since this would have enlarged too much the length of this paper. Of course, it would be straightforward to rebuild the full result for (\ref{osc-evol}), but this is far more complicated for (\ref{hom-P}) : we recall that there exists an extensive body of work dealing with the homogenization of the Dirichlet problem for fully nonlinear  first- and second-order partial differential equations in periodic, quasi-periodic, almost periodic and, more recently, stationary ergodic media. Listing references is beyond the scope of this paper but we tried to put several references in the bibliography (but our list is by no mean complete!).

Concerning results on oscillating boundary data, we were unable to find a lot of references~: in the papers of Alvarez and Bardi \cite{AB1,AB2} and Alvarez, Bardi and Marchi \cite{ABM}, the effects of oscillations in the initial data of parabolic equations are studied together with singular perturbations in the equations. For oscillating Dirichlet boundary conditions, the only related works are the ones of Gerard-Varet and Masmoudi \cite{GVM1,GVM2} : their results concern linear elliptic equations (or systems) in divergence form. By using energy methods, they are able to obtain results in polygonal domains but also the (difficult) case of smooth domains. It is worth pointing out that, despite of the very different methodological approaches, we face similar difficulties.

Our approach to recover the Dirichlet boundary condition consists, in both cases, in using a blow-up argument at any point of the boundary. For example, in the case of (\ref{hom-P}), if $x\in \dP$, we consider the functions $\ve$ defined by
$$ \ve (y):= \ue (x+\e y)\; ,$$
and we study their asymptotic behavior as $\e \to 0$. Under suitable assumptions on the $F_\e$'s, it is easy to see that, at least formally, their limit can be expressed through the solution $v$ of a Dirichlet problem in a half-space type domain, typically
\begin{equation}\label{Dpb-hstd}
\left\{
\begin{array}{rcl}
F (D^2 v, D v ,  y)  & = & 0  \quad \hbox{in  }P\, ,\\
v (y) & = & \psi (y)
\quad \hbox{on  }\dP \; ,
\end{array}
\right.
\end{equation}
for suitable functions $F$ and $\psi$. For example, for (\ref{hom-P}), if the restriction to $\dP$ of $g(y,x)$ is periodic, $\psi$ is just obtained from $g(y,x)$ by a translation in $y$ but this is more complicated if this condition is not satisfied (See the almost periodic case below).

In any case, Problem (\ref{Dpb-hstd}) plays a key role in recovering the Dirichlet boundary condition. The key fact is that -- under suitable conditions on $F$ (i.e. on the $\Fe$ or $\Ft$'s)--, (\ref{Dpb-hstd}) has a unique solution $v$ such that $v(y+Re)$ converges {\em uniformly for $y\in \bar P$} to a constant as $R \to +\infty$ and this constant provides the expected $\gb (x)$. Indeed, if the above heuristic argument is correct, then $\ue (x+\e Re)$ should be like $v(Re)$ and if this property holds locally uniformly in $x$, this means that $\ue (x)$ is like $\gb (x)$ on the set $\{y:\ y\cdot e=R\e \}$, and, knowing this behavior, it is easy to conclude.

This a priori simple program faces several difficulties : first the Dirichlet problem (\ref{Dpb-hstd}) is not always well-posed, even if we restrict ourselves to the (natural) framework of bounded solutions. To overcome this difficulty, we use a condition which already appear in the Neumann framework (cf. Barles, Da Lio, Lions, \& Souganidis \cite{BDLLS}).
Next the above convergence to a constant relies on a Liouville type property and requires additional properties on $F$. Finally we need a uniform convergence of $v(y+Re)$ as $R\to +\infty$ and except for the case when $\psi$ (and therefore $v$) has some periodicity property, this is a non trivial result.

The paper is organized as follows : in Section~\ref{dp-in-hstd}, we study Problem~(\ref{Dpb-hstd}) giving existence and uniqueness results together with the above mentioned property of the convergence to a constant. Section~\ref{appl:hom} describes the application of this result to Problem (\ref{hom-P}), while the last section provide the (analogue) approach to study Problem~(\ref{osc-evol}).

\section{The Dirichlet problem in a half-space type domain}\label{dp-in-hstd}

This section is devoted to the study of Problem (\ref{Dpb-hstd}).

Our first assumption is that $\psi$ is a continuous, periodic function defined on $\R$; this is a natural assumption since it is satisfied in all the applications we have in mind. As we mention in the introduction, the fact that the restriction of $\psi$ to $\dP$ is periodic in $y$ or not plays a role, at least at the technical level. In any case, the restriction of $\psi$ to $\dP$ is always almost periodic.

In the sequel, we will say that a real-valued function $\chi$, defined on either $P$ or $\dP$, is periodic with respect to a basis $\FF_{N-1} = (f_1, f_2,\cdots,f_{N-1})$ of $\dP$ if and only if
$$ \chi(y + z_1 f_1 + \cdots + z_{N-1}f_{N-1})= \chi(y)\;$$
for any $y \in P$ or $\dP$ and $z_1, \cdots, z_{N-1} \in \Z$, and we will write, in short, that $\chi$ is $\FF_{N-1}$-periodic. Throughout this section, we fix such a basis $\FF_{N-1} $ of $\dP$ and we complement it as a basis of $\R$ by adding a $N^{th}$ vector $f_N$. We define in an analogous way as above a $\FF_{N}$-periodic function in $\R$.\\

Concerning the nonlinearity $F$, we use the same assumptions as in \cite{BDLLS} without trying to state the most general conditions.

\medskip

\medskip
\leftline{(F0)\hfill
$\left\{\begin{array}{l}
\text{$F : \sym  \times\R \times\R \to\r$ is a continuous,  $\FF_{N}$-periodic function in $y$. }\\
\text{Moreover $F(M,p,y)$ is convex in $(M,p)$ and $F(0,0,y)= 0$ in $\R$.}\\
\end{array}\right.$}
\

\medskip
\leftline{(F1)\quad
$\left\{\begin{array}{l}
\text{$F$ is locally Lipschitz continuous on $\sym\times \R \times \R $
and there exists }\\
\noalign{\vskip6pt}
\text{a constant $K>0$ such that, for all $y_1,y_2\in \R$, $p,q \in \R$, and $M,N \in \sym$,}\\
\noalign{\vskip6pt}
 |F(M,p,y_1) - F(N,q,y_2)|\leq K
\{ |y_1-y_2|(1+|p|+|q| + ||M||+||N||) \\
\noalign{\vskip6pt}
\hskip2.2truein +|p-q| + ||M - N||\}\ ,
\end{array}\right.$ }

\medskip
\leftline{(F2)\quad
$\left\{\begin{array}{l}
\text{there exists $\kappa >0$ such that, for all $x\in \R$,
$p\in \R$,}\\
\noalign{\vskip6pt}
\text{and $M,N \in \sym $ with $N\geq 0$,}\\
\noalign{\vskip6pt}
F(M+N,p,x) - F(M,p,x) \leq -\kappa\tr(N) \; ,\\
\end{array}\right.$}
\

\medskip
The last assumption plays a key role in the uniqueness proof and we justify it below.

\medskip

\leftline{(F3)\quad
$\left\{\begin{array}{l}
\text{There exists a continuous, bounded from above subsolution $\wb$ of
}\\
\text{$F (D^2 \wb, D \wb, y) = 0$ in $P$, such that $\wb (y+Re) \to -\infty$ as $R \to + \infty$,}\\
\text{ uniformly for $y\in P$.}
\end{array}\right.$
}
\

\medskip
Our result is the

\begin{thm}\label{exist-uni} Assume (F0)-(F3). For any bounded, uniformly continuous function $\psi$, there exists a unique viscosity solution $v$ of (\ref{Dpb-hstd}). Moreover, if $\psi$ is  $\FF_{N-1}$-periodic, then $v$ is also $\FF_{N-1}$-periodic.
\end{thm}

\medskip
Before providing the proof of Theorem~\ref{exist-uni}, we comment Assumption (F3) which is the only non-classical hypothesis.

In fact, it is clear that (\ref{Dpb-hstd}) does not have a unique solution without any additional condition on $F$. Indeed one can look at the following simple example in dimension $1$
$$ -w'' -bw' = 0\quad\hbox{in  }(0,+\infty),\quad \hbox{$w(0)$ being given}\; ,$$
 where $b \in \r$. An immediate computation shows that $w$ is given by
$$ w(x)= w(0) + {\bar c}(\exp(-bx)-1)\; ,$$
where $ {\bar c}$ is any real constant. If $b > 0$, all these solutions are bounded and therefore we do not have uniqueness in this case. On the contrary, if $b<0$, the only bounded solution corresponds to $ {\bar c}=0$ and $b=0$ is a good case for uniqueness too. Therefore one needs to remove the $b>0$-case in order to have uniqueness and this is done through (F3) which is satisfied in the $b\leq 0$-case by $\wb (x):=- x$ and not satisfied for $b>0$.

How to interpret (F3)? In the linear case, our analysis uses a stochastic interpretation of pdes. We consider the case when the equation in Problem (\ref{Dpb-hstd}) writes
\begin{equation}
\label{lin}
- {\rm tr}[a(x)D^2v]-b(x)\cdot Dv = 0\quad\hbox{in  }P\; .
\end{equation}
where, for any $x\in \overline P$, $a(x)=\sigma(x)\sigma^T(x)$ for some $N \times p$-matrix $\sigma(x)$ and $p \le N$. If the drift $b$ and the diffusion matrix $\sigma$ satisfy classical assumptions, i.e. they are bounded and Lipschitz continuous functions on $\overline P$, it is well-known that the properties of Equation~(\ref{lin}) are connected with those of the solution $(X_t)_{t\geq 0} $ of the stochastic differential equation
\begin{equation}\label{ROe}
dX_t = b(X_t)dt +\sqrt2 \sigma (X_t)dB_t
\quad , \quad X_0=x \in P ,
\end{equation}
where $(B_t)_{t\geq 0}$ is an $N$-dimensional Brownian motion.

A simple idea in order to build the function $\wb$ would be to introduce the function
$$ W(x):= \Ex(\tau)\; ,$$
where $\tau$ is the first exit time of the trajectory $(X_t)_{t\geq 0} $ from $P$, i.e.
$$ \tau = \inf\{t\geq 0:\ X_t \notin P\}\; ,$$
and $ \Ex$ is the conditional expectation with respect to the event $\{X_0=x\}$.

When well-defined, the function $W$ satisfies
$$- {\rm tr}[a(x)D^2W]-b(x)\cdot DW = 1\quad\hbox{in  }P\; ,$$
$$ W = 0\quad\hbox{on  }\dP\; .$$
Therefore, a good candidate for $\wb$ in (F3) for (\ref{lin}) could be $-W$. Unfortunately, when (F3) is not satisfied, $W(x)\equiv +\infty$ (although $\tau$ is finite almost surely).

Indeed, a simple computation can be done in the above $1$-d example. We assume (F3) does not hold ($b>0$) and we prove $W=+\infty$. To do this, we compute the solution $W_R$ of
$$ -W_R'' -bW_R ' = 1\quad\hbox{in  }(0,R)\; ,$$
with $W_R (0)=W_R(R)=0$. In this case, we have replaced the first exit time from $(0,+\infty)$ by the exit time from the strip $(0,R)$ and the aim is, of course, to let $R$ tend to $+\infty$, since, obviously ,$W\geq W_R$ for all $R>0$. For $b > 0$, straightforward computations give
$$ W_R(x)= \frac1b\left(K_R(1-\exp(-bx))-x\right)\; ,$$
with, for large $R$, $K_R \sim R$. This proves that the expectation of the first exit time from the strip tends to infinity as $R$ tends to infinity.

Therefore (F3) means, in some sense, that the stochastic process $(X_t)_{t\geq 0} $ exits from $P$ ``rather quickly'' and this is a condition in order that the equation ``sees'' the boundary in a sufficient way to keep uniqueness.

An analogous difficulty, in the case of Neumann boundary condition, appears in \cite{BDLLS} where a similar assumption to (F3) is introduced. This assumption is the following

\medskip
\leftline{(F'3)\quad
$\left\{\begin{array}{c}
\text{there exists a bounded, Lipschitz continuous subsolution $\hat w$ of}\\
\noalign{\vskip6pt}
F (D^2\hat w, x,-e + D\hat w, y) = 0 \quad\hbox{in  }\R\; ,
\end{array}\right.$}
\

\medskip
(F'3) is stronger than (F3) since one can choose $\wb (y):=-e \cdot y + \hat w(y)$.

\medskip

Now we turn to the {\bf Proof of Theorem~\ref{exist-uni}.} The proof follows the classical steps of the Perron's Method introduced by Ishii \cite{I2} with the version ``up to the boundary'' by Da Lio \cite{DLCPDE}.

We first remark that $-||\psi||_\infty$ and $+||\psi||_\infty$ are respectively viscosity sub and supersolution of the (generalized) Dirichlet problem (\ref{Dpb-hstd}) because $F(0,0,x) \equiv 0$. Applying readily the method of \cite{DLCPDE}, we first build a discontinuous solution $v$ of this problem and then using assumptions (F1)-(F2), we deduce that, on $\dP$, no loss of boundary condition can occur and therefore
$$ v^* \leq \psi \quad\hbox{on  }\dP\quad \hbox{and} \quad v_* \geq \psi \quad\hbox{on  }\dP\; ,$$
where, as usual, $ v^*,  v_*$ denotes the upper and lower continuous envelope of $v$.

In order to deduce that $ v^*= v_*$ in $P$, a Strong Comparison Result is needed and this is where (F3) comes into play.

\begin{lemma}\label{comp} Assume (F0)-(F3). If $v_1$ is a bounded, usc viscosity subsolution of (\ref{Dpb-hstd}) associated to a bounded, uniformly continuous Dirichlet boundary condition $\psi_1$ and if $v_2$ is a bounded, lsc viscosity supersolution of (\ref{Dpb-hstd}) associated to a bounded, uniformly continuous  Dirichlet boundary condition $\psi_2$, then
$$||(v_1-v_2)^+||_\infty \leq ||(\psi_1-\psi_2)^+||_\infty \; .$$
\end{lemma}

The end of the proof is easy by using Lemma~\ref{comp} : first we use it with $v_1 = v^*$, $v_2 = v_*$ and $\psi_1=\psi_2=\psi$; this gives $ v^*\leq v_*$ in $P$ and therefore  $ v^*= v_*$ in $P$ since the opposite inequality is a direct consequence of the definitions of $ v^*, v_*$. Hence (\ref{Dpb-hstd}) has a continuous solution and uniqueness is again a straightforward consequence of Lemma~\ref{comp}.

Next the periodicity property comes from uniqueness since $v(y + z_1 f_1 + \cdots + z_{N-1}f_{N-1})$ and $v(y)$ are solutions with the boundary data $\psi (y + z_1 f_1 + \cdots + z_{N-1}f_{N-1})$ and $\psi(y)$ respectively; therefore they are equal if the former are.

For the {\bf proof of Lemma~\ref{comp}}, we first remark that we may assume without loss of generality that $\wb \leq -||\psi||_\infty$ on $\dP$; indeed $\wb$ is bounded from above and we can reduce to that case by subtracting a large enough constant. Then the proof consists in showing that, for any $0< \eta < 1$
$$ \eta v_1 + (1-\eta) \wb \leq v_2 \quad\hbox{in  }P\; .$$

The function $F$ being convex, the left-hand side of this inequality is a subsolution. Recalling the fact that $\wb (y) \to -\infty$ as $y\cdot e \to + \infty$, this inequality has to be proved only in a strip $\{y: 0\leq y\cdot e \leq R_\eta\}$.

The next argument is that, in the strip, we can build a strict subsolution for $F$: indeed, as a consequence of the uniform ellipticity and Lipschitz continuity in $p$ of $F$, it is enough to choose the function $y \mapsto \exp(Ly\cdot e)$, for $L>0$ large enough. It is worth pointing out that we may not even need the convexity at this stage : if $u$ is subsolution and if we set $u_\alpha (y):= u(y) + \alpha \exp(Ly\cdot e)$ for $0<\alpha \ll 1$, then, using (F1)-(F2), we have formally
\begin{align*}
F(D^2 u_\alpha (y), Du_\alpha (y), y)&=F(D^2 u(y) +\alpha L^2\exp(Ly\cdot e)e\otimes e, Du (y)+\alpha L\exp(Ly\cdot e)e, y),\\
& \leq F(D^2 u (y), Du (y), y)-\alpha \kappa L^2\exp(Ly\cdot e)+ K\alpha L\exp(Ly\cdot e),\\
& \leq -\alpha L\exp(Ly\cdot e) \bigl(\kappa L -K\bigr),
\end{align*}
which is clearly strictly negative (and even less than a strictly negative constant if $L>\kappa^{-1}K$.

At this point, the comparison becomes completely standard, applying the above argument to the subsolution $u=  \eta v_1 + (1-\eta)$ and using the arguments of Barles and Ramaswamy \cite{GBMR}. We conclude by letting first $\alpha$ tend to $0$ and then $\eta$ tend to $1$.$\sn$

\medskip
Now we turn to the behavior of the solution $v$ when $y\cdot e \to +\infty$. The  result is the
\begin{thm}\label{liouville} Assume (F0)-(F3). For any bounded, $\FF_{N-1}$-periodic continuous function $\psi$, if $v$ is the unique viscosity solution of (\ref{Dpb-hstd}), then $v(y+Re)$ converges uniformly for $y \in \bar P$ to a constant $\mu(\psi,F)$as $R \to +\infty$. The result is also true if $\psi$ is only an almost periodic continuous function provided that $F$ is independent of $y$.
\end{thm}

\noindent {\bf Proof.} We are going to prove the result in several steps.\\
{\bf Step 1 :} We show that the local uniform convergence holds, up to a subsequence, without assuming that $\psi$ is $\FF_{N-1}$-periodic.\\
We first notice that $v$ is uniformly bounded: more precisely, by Lemma~\ref{comp} (or by construction), $-||\psi||_\infty \leq v \leq ||\psi||_\infty$ on $\bar P$ since $-||\psi||_\infty$ and $||\psi||_\infty$ are respectively sub and supersolutions of (\ref{Dpb-hstd}). Next, by standard arguments (cf. Ishii and Lions \cite{il} or Barles \cite{japon}), $v_R (x):=v(x+Re)$ satisfies uniform $C^{0,\alpha}$-estimates. By Ascoli's Theorem and using the $\FF_N$-periodicity of $F$, one can extract from $(v_R)_R$ a subsequence which converge locally uniformly to a function $V$ which is a solution of
$$F (D^2 V, D V ,  x + \tau)   =  0  \quad \hbox{in  }\R\, ,$$
for some $\tau$ in $\R$.

It remains to prove that $V$ is constant, i.e. to obtain a Liouville type result for this equation. To do so, we use an idea of L. Rossi \cite{LR} ; the proof is done in two steps : in the first one, we show that $V$ is necessarely $\FF_{N}$-periodic and then we use the Strong Maximum Principle of Bardi and Da Lio \cite{BaDL} to conclude.

To perform the first step, we argue by contradiction : if $V$ is not $\FF_{N}$-periodic then there exists a vector $f_j$ and $x\in \R$ such that $V(x+f_j) \neq V(x)$. Without loss of generality, we can assume that $V(x+f_j) > V(x)$ (otherwise we can just change $f_j$ in $-f_j$) and we consider a sequence $x_k$ such that
$$V(x_k+f_j) - V(x_k) \to \sup_{\R}\,(V(x+f_j) - V(x))>0\; .$$
The sequence $(V(x_k+x))_k$ is also bounded and satisfies uniform $C^{0,\alpha}$-estimates, therefore we can extract a subsequence which converges locally uniformly, and the limit $\tilde V$ is a solution of the equation in $\R$ with in addition
$$\tilde V(f_j) - \tilde V(0) = \sup_{\R}\,(\tilde V(x+f_j) - \tilde V(x))>0\; .$$
Therefore the function $\tilde V(x+f_j) - \tilde V(x)$ achieves its maximum at $0$ but, by using the Strong Maximum Principle of Bardi and Da Lio \cite{BaDL}, this means that $\tilde V(x+f_j) - \tilde V(x)$ is constant and even equal to a positive constant $\delta$. This property is a contradiction with the fact that $\tilde V$ is bounded : indeed, for any $k \in \N$
$$ \tilde V(x+kf_j) = ( \tilde V(x+kf_j) - \tilde V(x+(k-1)f_j) )+ \cdots  ( \tilde V(x+f_j) - \tilde V(x)) + \tilde V(x)= k\delta + \tilde V(x)\; ,$$
and this is clearly impossible for large $k$ because $\tilde V$ satisfies the same $L^\infty$-bound as $v$.

Knowing that $V$ is $\FF_{N}$-periodic, this means that $V$ achieved its max (or min) on $\R$ at some point and again the Strong Maximum Principle of Bardi and Da Lio \cite{BaDL} allows to conclude that $V$ is constant.

\medskip
\noindent{\bf Step 2 : }The periodic case.\\
By Step 1, we know that there exists a sequence $R_k \to \infty$ such that $v(y+R_k e)$ converges locally uniformly to a constant but, by taking $y \in \dP$ and using the periodicity, we see that it is uniform on the hyperplan $\dP_k:=\{x \cdot e = R_k \}$. Hence $y\mapsto v(y+R_k e)$ is a solution of a problem which is similar to (\ref{Dpb-hstd}), $F$ being changed in $F(\cdot,\cdot,x+R_k e)$ and is uniformly close to a constant $\lambda$ on $\dP_k$. Moreover the constant $\lambda$ solves $F=0$ by (F0) and therefore, by Lemma~\ref{comp}
$$||v(y+R_k e)-\lambda||_{L^\infty(P)} \leq ||v(y+R_k e)-\lambda||_{L^\infty(\dP)} \; .$$
This inequality obviously shows that $v(y+R_k e)$ converges to $\lambda$ uniformly on $P$, and $v(y+R e)$ too.

\medskip
\noindent{\bf Step 3 : }The almost periodic case.\\
We first recall that $\psi$ being almost periodic, for any $\e>0$, there exists a relatively dense set $T_\e$ of $\dP$, such that
$$||\psi (x+\tau) - \psi(x)||_{L^\infty(\dP)} \leq \e\quad \hbox{for any  }\tau \in T_\e\; ,$$
and we also recall that ``relatively dense set'' means that there exists $S(\e)$ such that
$$ \dP = \bigcup_{\tau \in T_\e}\,B(\tau,S(\e))\; .$$
We define $R_k$ as in the previous step. Since in this case $F$ is independent of $y$, for any $\tau \in T_\e$, $v(\cdot+\tau)$ and  $v$ are solutions of the $F$-equation and using again Lemma~\ref{comp}, we deduce
$$||v(y+\tau)- v(y)||_{L^\infty(P)} \leq ||\psi (y+\tau) - \psi(y)||_{L^\infty(\dP)} \leq \e  \; .$$

In particular, this inequality implies
$$||v(y+\tau + R_k e)- v(y+ R_k e)||_{L^\infty(\dP)} \leq ||\psi (y+\tau) - \psi(y)||_{L^\infty(\dP)} \leq \e  \; .$$
If, by the local uniform convergence, $v(y+ R_k e) \to \lambda$, for some constant $\lambda$, on the ball $B(0,S(\e))$, we deduce from the previous inequality that
$$  ||v(y+\tau + R_k e)- \lambda||_{L^\infty(\dP)} \leq \e + o_k(1) \; ,$$
where $o_k(1) \to 0$ as $k\to +\infty$ and therefore, recalling the property of the set $T_\e$ and using the same argument as in Step 2
$$  ||v(y + R_k e)- \lambda||_{L^\infty(P)} \leq \e + o_k(1) \; .$$
Playing with the parameters $k$ and $\e$, one easily proves the global convergence property and the conclusion follows.
$\sn$

\begin{remark}\label{goodvar} There exists several variants of the results of this section~: in particular, we point out the case when $F(M,p,x)$ does not depend on $p$ and is therefore only a function of $M$ and $x$. In this case, (F3) is automatically satisfied by $ \bar w (y)=-y\cdot p$ and, as a consequence, (F0) can be reduced to continuity and $F(0,x)=0$, i.e. the convexity can be dropped.\end{remark}

\section{Application to Homogenization Problems in half-space type domains}\label{appl:hom}

We consider the homogenization problem (\ref{hom-P}). Under suitable assumptions, it is expected that the $\ue$'s converge to the solution $u$ of a Dirichlet problem of the form
\begin{equation}\label{hom-P-lim}
\left\{
\begin{array}{rcl}
{\overline F} (D^2 u, D u , u ,  x)  & = & 0  \quad \hbox{in  }P\, ,\\
\displaystyle u  & = & \gb (x)
\quad \hbox{on  }\dP \; ,
\end{array}
\right.
\end{equation}
for some homogenized, uniformly elliptic, nonlinearity ${\overline F}$ and boundary condition $\gb$.
We recall that our aim is just to identify the Dirichlet boundary condition $\gb$ and therefore we are going to concentrate on this issue. We refer to \cite{AB1, AB2, BDLLS, E1, E2} and references therein for the other ingredients allowing to treat the full problem.

As mentioned in the introduction, in order to recover the boundary condition at $x\in \dP$, we consider the function $ \ve (y):= \ue (x+\e y)$ which solves
$$ F_\e (\e^{-2} D^2 \ve, \e^{-1} D \ve , \ve , \xoe +y, x+\e y)   =  0  \quad \hbox{in  } P\, ,$$
with
$$ \ve (y):= g (\xoe +y, x+\e y) \quad \hbox{on  }\dP \, .$$
The aim is to show that $\ve$ converges locally uniformly to a solution of a problem like (\ref{Dpb-hstd}) and to do so, we use the following assumptions.

\medskip
\leftline{(F4)\;
$\left\{\begin{array}{l}
\text{There exists $k>0$ such that $\e^k F_\e (\e^{-2} M, \e^{-1} p , u , y, x+\e y)$ converges} \\
\noalign{\vskip6pt}
\text{locally uniformly in $(M,p,u,y,x)$ to a function $F(M,p,y,x)$ satisfying }\\
\noalign{\vskip6pt}
\text { (F0), (F1), (F2), (F3) as a function of $M,p,y$.}\\
\end{array}\right.$ }
\ \\
\medskip
\leftline{(F5)\quad
$\left\{\begin{array}{l}
\text{There exists $k>0$ such that $\e^k F_\e (\e^{-2} M, \e^{-1} p , u , y, x+\e y)$ converges locally}
\\
\noalign{\vskip6pt}
\text{ uniformly in $(M,p,u)$ and uniformly in $(x,y)$ to a function $F(M,p,x)$ }\\
\noalign{\vskip6pt}
\text {independent of $y$ satisfying (F0), (F1), (F2), (F3) as a function of $M,p,y$.}
\end{array}\right.$ }

Our result is the

\begin{thm}\label{bc-Dir-lim} Assume that there exists $\e_0 >0$ such that, for any $0< \e \leq \e_0$, there exists a solution $\ue$ of (\ref{hom-P}) and that the $\ue$ are uniformly bounded. At $x\in \dP$, we have $$\gb (x)=\mu(g(\cdot,x),F(\cdot,\cdot,x))\; ,$$ where $\mu$ is defined in Theorem~\ref{liouville} in both following cases
\\
{\bf (i)} If $g(\cdot,x)$ is a continuous, $\FF_{N-1}$-periodic function and (F4) holds.\\
{\bf (ii)} If $g(\cdot,x)$ is a bounded continuous, almost-periodic function and (F5) holds.
\end{thm}

\noindent {\bf Proof.} It is based on the following lemma.
\begin{lemma}\label{lemma-inter}Under the assumptions of Theorem~\ref{bc-Dir-lim}, for any $\delta >0$, there exists $R_0$ such that, for $R\geq R_0$ and for any sequence $\xe \to x$
$$ \limsup_\e |\ue (\xe+\e R e)-\gb(x)| \leq \delta\; .$$
 \end{lemma}

Using this lemma, the passage to the limit in the boundary condition at $x$ is rather easy because it is enough to consider the functions $x \mapsto \tilde u^\e (x):=\ue (x+\e R e)$ which satisfy Equation~(\ref{hom-P}) in $P$ (with $\xoe$ changed into $\xoe+R e$ and $x$ in $x+\e R e$) and, on the boundary of $P$, we have, for $x'$ close enough to $x$,
$\gb(x)- \delta \leq \tilde u^\e (x') \leq \gb(x) + \delta \; .$

 The passage to the limit follows from standard arguments, letting finally $\delta$ tend to $0$.\\

 It is worth pointing out that the half-relaxed limits of $\ue$ and $\tilde u^\e$ are the same in $P$ but not on $\partial P$ : the translation by $\e R e$ in the definition of $\tilde u^\e$ is used to get rid of the boundary layer for $\ue$ and to recover the proper boundary condition. We also note that the uniform ellipticity of all the equations we consider implies that the Dirichlet boundary conditions are satisfied in the usual sense (and not only in the viscosity solutions' sense).\sn \\

\noindent{\bf Proof of Lemma~\ref{lemma-inter}.} For the proof in the case {\bf (i)}, we denote by $\Pi$ the set of points of the form $x_1 f_1 + \cdots + x_{N-1}f_{N-1}$ with $0\leq x_i\leq 1$ for $i=1,\cdots,N-1$.

Since $x \in \dP$, $ \xoe \in \dP$ and the functions $F_\e,g$ being $\FF_{N-1}$-periodic in $y$, we can assume without loss of generality that the $\xoe$ are in $\Pi$. Moreover, since $\Pi$ is compact, we can even assume that $\xoe$ is converging to some $\tau \in \Pi$.

Applying the half-relaxed limit method to pass to the limit in the $\ve$-problem (in fact to the subsequence along which $\xoe$ is converging to $\tau$), we easily prove that, up to a subsequence, $\ve$ converges locally to the unique solution of (\ref{Dpb-hstd}) associated with the non-linearity $F(M,p,\tau + y,x)$ and the Dirichlet boundary condition $g(\tau +y,x)$. We point out that the subsequence just selects the limiting $\tau$.

% The only difficulty (and this is a key one since this prevents us to provide results in general domains) comes from the dependence in $\xoe$ in the boundary condition $g (\xoe +y, x+\e y)$. Of course, the $\e y$-term plays no role and we can look at $g (\xoe +y, x)$ and we see this function as a $\Z^N$-periodic function of $y$ in the whole space $\R$. Therefore we can extract uniformly converging subsequences ; we denote by $\tilde g(y,x)$ one of the possible limits.
%
%Coming back on $\dP$, the important point is that the solutions of (\ref{Dpb-hstd}) associated respectively to $g (\xoe +y, x)$ are $v (\xoe +y, x)$ where $v (y, x)$ is the solution associated to $g(y,x)$ since $F$ does not depend on $y$ and therefore all these functions have the same behavior when $y\cdot e \to +\infty$, i.e. the same $\mu(g(\cdot,x),F(\cdot,\cdot,x))$ and so has $\tilde g(y,x)$ by uniform convergence.

But, if $v$ is the solution of this problem with $\tau = 0$, then, by uniqueness, the solution for a general $\tau$ is just $v(\tau + \cdot)$ and all these functions have the same behavior when $y\cdot e \to +\infty$, namely they converge uniformly to $\mu(g(\cdot,x),F(\cdot,\cdot,\cdot,x))$.

The result is therefore just a consequence of Theorem~\ref{liouville} together with the properties recalled in the previous paragraph.

Finally the local uniform convergence comes the fact that, in the above argument, we can replace $\ve$ by $\tilde v_\e (y):= \ue (\xe+\e y)$ for any sequence $(\xe)_\e$ converging to $x$.

For the proof of {\bf (ii)}, we almost argue in the same way : because of the assumptions on $g$, up to a subsequence, $g (\xoe +y, x)$ converges uniformly on $\dP$ to some function $\tilde g(y,x)$ which is almost periodic in $y$. As a consequence, for the same subsequence, $g (\xoe +y, x + \e y)$ is converging locally uniformly to $\tilde g(y,x)$ and applying the half-relaxed limit method to pass to the limit in the $\ve$-problem (in fact to the subsequence), we obtain the local uniform convergence to
the unique solution $v$ of (\ref{Dpb-hstd}) associated with the non-linearity $F(M,p,x)$ and the Dirichlet boundary condition $\tilde g(y,x)$.

In order to conclude, we need to show that the behavior when $y\cdot e \to +\infty$ of the functions $v$ is independent of the subsequence. To do so, we have just to consider the unique solution $w$ of (\ref{Dpb-hstd}) associated with the non-linearity $F(M,p,x)$ and the Dirichlet boundary condition $g(y,x)$. Since $F(M,p,x)$ is independent of $y$, $w (\xoe +y)$ is the unique solution of (\ref{Dpb-hstd}) associated to $g (\xoe +y, x)$ and therefore $w (\xoe +y)$ is converging locally uniformly to $v$. But, by Theorem~\ref{liouville}, the functions $w (\xoe +y)$ have the same behavior when $y\cdot e \to +\infty$ and so has $v$. And the conclusion follows as in case {\bf (i)}.
$\sn$

\section{Examples and Remarks}

In order to illustrate our assumptions and to emphasize the difficulty, we consider two linear problems. The first one is
\begin{equation}\label{eq1.1}
- \tr (A(x,\xoe)D^2 \ue) - b(x,\xoe)\cdot D \ue+ \ue  =  f(x, \xoe)  \quad \hbox{in  }P\, ,
\end{equation}
while the second one, which is more singular, is
\begin{equation}\label{eq1.2}
 - \tr (A(x,\xoe)D^2 \ue) - \displaystyle {\e}^{-1}b(x,\xoe)\cdot D \ue +  \ue
  =  f(x, \xoe)   \quad \hbox{in  }P\, ,
\end{equation}
 which has, as
a particular case,  the divergence form equation
$$ - {\rm div}( A(x,\xoe) D\ue) + \ue = f(x,\xoe) \quad \hbox{in  } P\; .$$

In order to simplify the exposure, we assume that the functions $A(x,y), b(x,y), f(x,y)$ are bounded, Lipschitz continuous, $\FF_{N}$-periodic functions (we can slightly weaken these regularity assumptions)
and there exist $0< \nu \leq \Lambda$ such that, for any $x\in P,y \in \R$
$$ \nu Id_N \leq A(x,y)\leq \Lambda Id_N\; .$$

These equations are, of course, associated with the Dirichlet boundary condition of (\ref{hom-P}).

Concerning (\ref{eq1.1}), we have to check (F4) and to compute $F$ : choosing $k=2$ which is in general the right choice, we obtain
$$
F(M,p,y,x):= - \tr (A(x,y)M) \, ,
$$
therefore independent of $p$. Taking into account Remark~\ref{goodvar}, it is clear that (F0)-(F3) holds and with suitable assumption on $g$, Theorem~\ref{bc-Dir-lim} {\bf (i)} applies readily or {\bf (ii)} if $A$ is independent of $y$.

Concerning (\ref{eq1.2}), with the same choice of $k$, the computation of $F$ gives
$$
F(M,p,y,x):= - \tr (A(x,y)M) - b(x,y)\cdot p \, .
$$
The checking of (F0)-(F2) remains easy but, as we explain it in Section~\ref{dp-in-hstd}, it is not obvious to give general conditions ensuring that (F3) holds. We recall that the above divergence form equation enters into this framework and its treatment differs from (for example) the simplest case where we assume
$$ b(x,y)\cdot e \leq 0\quad \hbox{for any  }x\in P,\ y \in \R,$$
where (F3) is satisfied by $\bar w (y) = -e\cdot y$.

If we assume that (F3) holds and with suitable assumption on $g$, Theorem~\ref{bc-Dir-lim} {\bf (i)} applies readily or {\bf (ii)} if $A,b$ are independent of $y$.

Now we turn to a key question : why not giving such results in smooth bounded domains? If $\OO$ is a smooth bounded domain and if we want to prove the analogue of Theorem~\ref{bc-Dir-lim}, replacing $P$ by $\OO$, the blow up argument near a point $x \in \OO$ leads us to consider Problem (\ref{Dpb-hstd}) in $P:=\{y:\ y\cdot n(x) \leq 0 \}$ where $n(x)$ stands for the unit, normal to $\partial \OO$ pointing outward to $\OO$. Then we have two cases : (i) either $n(x)$ has a rational slope and the restriction of the, say, $\Z^N$-periodic in $y$ function $g(y,x)$ to $P$ is indeed $\FF_{N-1}$-periodic for a suitable choice of vectors $f_i$ ($1 \leq i \leq N-1$) (ii) or $n(x)$ has no rational slope then the restriction of $g(y,x)$ to $P$ is only almost periodic and we are stuck with rather restrictive assumptions on the equation, typically (\ref{eq1.1}) with $A$ being independent of $y$.

The above examples and remarks show that the results of \cite{GVM1} solves all the main difficulties of this type of problems.

\section{The parabolic case : problems with oscillations in time}

We consider the initial-boundary value problem (\ref{osc-evol}). Here again the aim is just to recover the limiting Dirichlet boundary condition on $\dO$, namely $\overline g$.

To do so, for $(x,t)\in \dO \times (0,T)$, we introduce the function
$$\ve (y,s):= \ue (x + \e y, t + \e^2 s)\; .$$
The function $\ve$ satisfies the equation
$$\ve_s + \e^2 \Ft (\e^{-2} D^2 \ve , \e^{-1} D \ve, x + \e y, t + \e^2 s)  =  0
\quad \hbox{in  }\e^{-1} (\O-x) \times (-\e^{-2}t ,\e^{-2} (T-t)) \, ,$$
with $g(x + \e y, \e^{-2}t + s)$ as Dirichlet boundary condition.\\

To study the asymptotic behavior of $\ve$, we first assume

\medskip
\leftline{(F6)\quad
$\left\{\begin{array}{l}
\text{The functions $\e^2 \Ft (\e^{-2} M, \e^{-1} p , x ,t)$ converge locally uniformly to}
\\
\noalign{\vskip6pt}
\text{  a function $\bar F (M,x,t)$. Moreover $M\mapsto \bar F (M,x,t)$ satisfies (F1), (F2)}
\end{array}\right.$ }

\medskip

Using this assumption and ignoring, for the moment, the term  $\e^{-2}t$, a formal passage to the limit gives the analogue of problem (\ref{Dpb-hstd}), namely
\begin{equation}\label{Dpb-hstd-evol}
\left\{
\begin{array}{rcl}
v_s + \bar F (D^2 v, x , t)  & = & 0  \qquad \hbox{in  }P\times \r\, ,\\
v  & = & g (x,s)
\quad \hbox{on  }\dP \times \r\; ,
\end{array}
\right.
\end{equation}
where $P := \{y:\ y\cdot n(x)\leq 0\}$, $n(x)$ being the unit outward normal vector to $\dO$ at $x$. It is worth pointing out that it is actually a problem in the variables $(y,s)$, $(x,t)$ being considered as parameters.

This analysis leads to consider the following general problem which is the analogue of  (\ref{Dpb-hstd})
\begin{equation}\label{Dpb-hstd-model}
\left\{
\begin{array}{rcl}
v_s + \tilde F (D^2 v)  & = & 0  \quad \hbox{in  }P\times \r\, ,\\
v  & = & \varphi (s)
\quad \hbox{on  }\dP \times \r\; ,
\end{array}
\right.
\end{equation}
where $P$ is the half-space $\{y:\ y\cdot e \geq 0\}$ and $\varphi$ is a continuous, $1$-periodic function in $\r$.

The analogue of Theorems~\ref{exist-uni} and \ref{liouville} is the

\begin{thm}\label{exist-uni-liouv} For any continuous, $1$-periodic function $\varphi$, there exists a unique viscosity solution $v$ of (\ref{Dpb-hstd-model}). This function is $1$-periodic in $s$ and depends on $y$ through $y\cdot e$ only. Finally, $v(y\cdot e , s)$ converges uniformly to a constant $\mu(\varphi,\tilde F)$ as $y\cdot e \to + \infty$.
\end{thm}

\noindent{\bf Proof.} We just sketch it since it is an easy adaptation of the proofs of Theorems~\ref{exist-uni} and \ref{liouville}. Here are the key arguments.\\
(i) The comparison of a (bounded) subsolution $v_1$ and a (bounded) supersolution $v_2$ can be done the following way : since $\tilde F$ does not depend on $Dv$, the function $(y,s) \mapsto v_1 (y,s)-\alpha y\cdot e$ is still a subsolution of (\ref{Dpb-hstd-model}), for any $\alpha >0$ and we have to show that $v_1 (y,s)-\alpha y\cdot e \leq v_2 (y,s)$ in $P\times \r$. This replaces the convexity arguments of the proof of Theorem~\ref{exist-uni} and, as in the proof of Lemma~\ref{comp}, we just have to prove that this inequality holds in the strip $\{y:\ 0\leq y\cdot e \leq R(\alpha)\}$, for some large enough $R(\alpha)$. The other arguments apply readily, in particular the change with the function $y \mapsto \exp(-Ly\cdot e)$, for $L>0$ large enough, which uses only the ellipticity in the space variable.\\
(ii) The facts that $v$ is $1$-periodic in $s$ and depends on $y$ only through $y\cdot e$ are immediate consequences of the uniqueness and the fact that $\tilde F$ does not depend either on $y$ or on $s$; this allows us to use translations with respect to these variables, in directions which are orthogonal to $e$.\\
(iii) The Liouville type result is easy to obtain because $\tilde F$ only depends on $D^2v$ : either we repeat the arguments of the proof of Theorem~\ref{liouville} or we can show the result through local $C^{0,\alpha}$ bounds : if $R = y\cdot e$ is large enough, it is easy to see that the local $C^{0,\alpha}$-norms in space in $B(y,R/2)$ tend to $0$ as $R\to +\infty$, uniformly in time, and so does the modulus of continuity in time by an argument of Barles, Biton and Ley \cite{BBL} .$\sn$

\begin{thm} \label{bc-Dir-lim-evol}
Assume (F6) and that (\ref{osc-evol}) have solutions which are uniformly bounded for $\varepsilon$ small enough. Then the boundary condition for the limiting problem is given, at $(x,t)\in \partial O\times (0,T)$, by $\bar g (x,t)=\mu(g(x,\cdot),\tilde F (\cdot,x,t))$ where $\mu$ is defined in Theorem~\ref{exist-uni-liouv}.
\end{thm}

We conclude our article by the analogue of Lemma~\ref{lemma-inter} from which Theorem~\ref{bc-Dir-lim-evol} can be easily deduced.

\begin{lemma} Assume (F6). For any $\delta >0$, there exists $R_0$ such that, for $R\geq R_0$ and for any sequences $(\xe,\te) \to (x,t)$
$$ \limsup_\e |\ue (\xe - \e R n(x), \te)-\bar g (x,t) | \leq \delta\; .$$
\end{lemma}

We skip the proof since it is analogous and even easier than the proof of Lemma~\ref{lemma-inter}: the periodicity in time allows to treat the bad $\e^{-2}t$-term by localizing it in a periodic cell.


\begin{thebibliography}{99}

\bibitem{AB1}
Alvarez, O. and  Bardi, M.:
{\sl~Viscosity solutions methods for singular perturbations in deterministic
and stochastic control,}
SIAM J. Control Optim. 40 (2001/02)  1159--1188.

\bibitem{AB2}  Alvarez, O. and  Bardi, M.:
{\sl~Singular perturbations of nonlinear degenerate parabolic PDEs: a
general convergence result,}
Arch. Ration. Mech. Anal. 170 (2003)  17--61.

\bibitem{ABM}  Alvarez, O., Bardi, M., Marchi, C. :
Multiscale problems and homogenization for second-order Hamilton-Jacobi equations.
J. Differential Equations 243 (2007), no. 2, 349Ð387.

\bibitem{A3} Arisawa, M.:
{\sl~Long time averaged reflection force and homogenization of oscillating
Neumann boundary conditions,}
Ann. Inst. H. Poincar\'e Anal. Lin\'eaire 20 (2003)  293--332.

\bibitem{AL} Arisawa, M. and  Lions, P.-L.:
{\sl~On ergodic stochastic control,  }
Comm. Partial Differential Equations 23 (1998)  2187--2217.

\bibitem{BaDL} Bardi, M. \& Da Lio, F.: {\sl On the strong maximum principle for fully nonlinear degenerate elliptic equations,}  Arch. Math. (Basel) 73 (1999), no. 4, 276-285.


\bibitem{bb} Barles, G.:
{\sc Solutions de viscosit\'e des \'equations de Hamilton-Jacobi}.
Collection ``Math\'ematiques et Applications'' de
la SMAI, n$^\circ$17, Springer-Verlag (1994).

\bibitem{b2}
Barles, G.:
{\sl Nonlinear Neumann boundary conditions for
quasilinear, degenerate elliptic equations and applications,}
Journal of Diff. Eqs. { 154} (1999) 191--224.

\bibitem{japon}
Barles, G.:  $C^{0,\alpha}$-regularity and estimates for solutions of elliptic and
parabolic equations by the Ishii \& Lions method.

\bibitem{BBL} Barles, G., Biton S. and Ley O. :
A geometrical approach to the study of unbounded solutions of
quasilinear parabolic equations, Arch.  Rational Mech.  Anal.
{\bf 162} (2002), 287--325.

\bibitem{BDL-eb} Barles, G. and  Da Lio, F.:
{ \sl On the boundary ergodic problem for fully nonlinear equations
in bounded domains with general nonlinear Neumann boundary conditions},
Ann. IHP, Analyse non lin\'eaire, to appear.

\bibitem{BDL}
Barles, G.   and   Da Lio, F.:
{\sl Local $C^{0,\alpha}$ estimates for viscosity solutions of
Neumann-type boundary value problems}, preprint.

\bibitem{BDLLS}
Barles, G., Da Lio, F., Lions, P.-L., Souganidis, P. E.,
{\sl Ergodic problems and periodic homogenization for fully nonlinear equations in half-space type domains with Neumann boundary conditions},
{Indiana Univ. Math. J.},
{\bf 57}, (2008), N. 5, 2355 -- 2376.

\bibitem{BaLi}  Barles, G. and  Lions, P.-L.:
{\sl Remarques sur les probl\`emes de reflexion obliques,}
C. R. Acad. Sci. Paris, t. 320 Serie I (1995)  69--74.

\bibitem{GBMR}  Barles, G. and  Ramaswamy, M.:
{\sl Sufficient structure conditions for uniqueness of viscosity solutions
of  semilinear and quasilinear equations,}
NoDEA, to appear.

\bibitem{BS1} Barles, G. and Souganidis, P.E.:
{\sl On the large time behaviour of solutions of Hamilton-Jacobi equations,}
SIAM J. Math.  Anal.  31 (2000) 925--939.

\bibitem{BS2} Barles, G. and Souganidis, P.E.:
{\sl Space-time periodic solutions and long-time behavior of solutions
to quasi-linear parabolic equations,}
SIAM J. Math. Anal. 32 (2001) 1311--1323.

\bibitem{BLP}
Bensoussan, A., Lions, J.-L., and Papanicolaou, G.:
{\tt Cited on page 3.}

\bibitem{B} Bensoussan, A.:
{\sc Perturbation methods in optimal control,}
Translated from the French by C. Tomson. Wiley/Gauthier-Villars Series
in Modern Applied Mathematics.
John Wiley \& Sons, Ltd., Chichester; Gauthier-Villars, Montrouge, 1988.

\bibitem{CSW}
Caffarelli, L.A., Souganidis, P.E. and Wang, L.:

\bibitem{cil}
Crandall M.G., Ishii, H. and  Lions, P.-L.:
{\sl User's guide to viscosity solutions of second order Partial differential
equations}, Bull. Amer. Soc. {27} (1992)  1--67.

\bibitem{DLCPDE} F. Da Lio: {\sl Strong Comparison Results for
Quasilinear Equations in Annular Domains and Applications},  Comm. in
PDE
{\bf 27}  (1\& 2), 283-323 (2002).


\bibitem{E1} Evans, L.C.:
{\sl The perturbed test function method for viscosity solutions of
nonlinear PDE,}
Proc. Roy. Soc. Edinburgh Sect. A 111 (1989)  359--375.

\bibitem{E2} Evans, L.C.:
{\sl Periodic homogenisation of certain fully nonlinear partial differential
equations,}
Proc. Roy. Soc. Edinburgh Sect. A 120 (1992)  245--265.

\bibitem{GVM1} Gerard-Varet, D. and Masmoudi, N. : {\sl Homogenization and Boundary Layer}. Preprint.

\bibitem{GVM2} Gerard-Varet, D. and Masmoudi, N. :
{\sl Homogenization in polygonal domains}.  J. Eur. Math. Soc. (JEMS) 13 (2011), no. 5, 1477Ð1503.

\bibitem{I1} Ishii, H.:
{\sl Almost periodic homogenization of Hamilton-Jacobi equations,}
International Conference on Differential Equations,
Vol. 1, 2 (Berlin, 1999), 600-605, World Sci. Publishing, River Edge, NJ, 2000.

\bibitem{I2} Ishii, H.:
{\sl Perron's method for Hamilton-Jacobi Equations,}
Duke  Math. J. {55} (1987) 369--384.

\bibitem{I4}  Ishii, H.:
{\sl Fully nonlinear oblique derivative
problems for nonlinear second-order elliptic PDE's,}
Duke Math.  J. { 62} (1991) 663--691.

\bibitem{il} Ishii, H. and Lions, P.-L.:
{\sl  Viscosity solutions of fully nonlinear second-order
elliptic partial differential equations,}
{J. Differ. Equations} {83}  (1990) 26--78.

\bibitem{LPV} Lions, P.-L., Papanicolaou, G., Varadhan, S.R.S.:
unpublished preprint.

\bibitem{LS1} Lions P.-L and Souganidis, P.E.: Homogenization of degenerate second-order PDE in periodic and almost periodic environments and applications. Ann. Inst. H. PoincarŽ Anal. Non LinŽaire 22 (2005), no. 5, 667Ð677.


\bibitem{LS2} Lions P.-L and Souganidis, P.E.:
Homogenization of "viscous'' Hamilton-Jacobi equations in stationary ergodic media.
Comm. Partial Differential Equations 30 (2005), no. 1-3, 335Ð375.

\bibitem{LS3} Lions P.-L and Souganidis, P.E.: Stochastic homogenization of Hamilton-Jacobi and "viscous''-Hamilton-Jacobi equations with convex nonlinearitiesÑrevisited. Commun. Math. Sci. 8 (2010), no. 2, 627Ð637.

\bibitem{LiSz} Lions P.-L and Sznitman, A.S.:
{\sl Stochastic differential equations with
reflecting boundary conditions, }
Comm. Pure and Applied Math.  37	
(1984) 511--537.

\bibitem{LR} Rossi, L. : personal communication.

\bibitem{HT} Tanaka, H.:
Homogenization of diffusion processes with boundary conditions.
Stochastic analysis and applications, 411--437, Adv. Probab.
Related Topics 7 Dekker, New York, 1984.

\end{thebibliography}
\end{document}